\documentclass[11pt, leqno]{article}
\usepackage{amsfonts}
\usepackage{amsmath}
\usepackage{amsthm}
\usepackage{amssymb}
\def\b{\mathbb }

\def\phi{\varphi }

\swapnumbers
\theoremstyle{plain}
\newtheorem{theorem}{Theorem}[section]
\newtheorem{corollary}[theorem]{Corollary}
\newtheorem{lemma}[theorem]{Lemma}
\newtheorem{proposition}[theorem]{Proposition}
\theoremstyle{definition}
\newtheorem{definition}[theorem]{Definition}
\theoremstyle{remark}
\newtheorem*{remark}{Remark}

\newtheorem{examples}[theorem]{Examples}

\numberwithin{equation}{section}

\begin{document}
\title{Generalized Hermite Polynomials and the Heat Equation for Dunkl 
Operators}
\author{Margit R\"osler\thanks{This paper was written  while the  author 
held a Forschungsstipendium of the DFG at the University of Virginia, 
Charlottesville, USA.} \\
 Mathematisches Institut, Technische Universit\"at M\"unchen\\
 Arcisstr. 21, D-80333 M\"unchen, Germany\\
 Current address:\\
 Department of Mathematics, University of Virginia, Kerchof Hall\\
 Charlottesville, VA 22903, USA\\
e-mail: roesler@mathematik.tu-muenchen.de\\
}

\date{}

\maketitle

\begin{abstract} Based on the theory of Dunkl operators, this paper presents 
a general concept of multivariable Hermite polynomials and Hermite functions 
which are associated with finite reflection groups on $\b R^N$.  The 
definition and properties of these generalized Hermite systems extend 
naturally those of their classical counterparts; partial derivatives and the 
usual exponential kernel are here replaced by Dunkl operators and the 
generalized exponential kernel $K$ of the Dunkl transform. In case of the 
symmetric group $S_N$, our setting includes the polynomial eigenfunctions of 
certain Calogero-Sutherland type operators. The second part of this paper is 
devoted to the heat equation associated with Dunkl's Laplacian. As in the 
classical case, the corresponding Cauchy problem is governed by a positive 
one-parameter semigroup; this is assured  by a maximum principle for the 
generalized Laplacian. The explicit solution to the Cauchy problem involves 
again the  kernel $K,$ which is, on the way, proven to be nonnegative for 
real arguments. 
\end{abstract}



\smallskip

\noindent
 Running title: Hermite polynomials and the heat equation for Dunkl operators.

\section{Introduction}

Dunkl operators are differential-difference operators associated with a 
finite reflection group, acting on some Euclidean space. They provide a 
useful framework for the study of multivariable analytic structures which 
reveal certain reflection symmetries. 
  During the last years, these operators have gained considerable interest in 
various fields of mathematics and also in physical applications; they are, 
for example,  naturally connected with certain Schr\"odinger operators for 
Calogero-Sutherland-type quantum many body systems, see \cite{LV} and 
\cite{BF2}, \cite{BF3}. 
 
For a finite reflection group  $G\subset O(N,\b R)$  on $\b R^N$  the 
associated Dunkl operators are defined as follows:  For $\alpha \in \b R^N$, 
denote by $\sigma_\alpha$ the reflection corresponding to $\alpha$, i.e. in 
the hyperplane orthogonal to $\alpha$.  It is given by
\[ \sigma_\alpha (x) = x- 2\frac{\langle\alpha,x\rangle}{|\alpha|^2}\,
\alpha,\] 
where  $\langle.,.\rangle$ is the Euclidean scalar product on $\b R^N$ and 
$|x|:= \sqrt{\langle x,x\rangle}$. (We use the same notations for the 
standard Hermitian inner product and norm on $\b C^N$.) Let $R$ be the  root 
system 
associated with the reflections of $G$, normalized such that 
$\langle\alpha,\alpha\rangle = 2$ for all $\alpha\in R$. 
Now choose a multiplicity function $k$ on the root system $R$, that is, a 
$G$-invariant function $k:R\to \b C$, and fix some positive subsystem $R_+$ 
of $R$.  The Dunkl operators $T_i \>(i=1,\ldots,N)$ on $\b R^N$ associated 
with $G$ and $k$ are then given by 
\[ T_if(x) := \partial_if(x) + \sum_{\alpha\in R_+} k(\alpha)\,\alpha_i\cdot 
\frac{f(x)-f(\sigma_\alpha x)}{\langle\alpha,x\rangle}\,,\>\>f\in 
C^1(\b R^N);\]
here $\partial_i$ denotes the $i$-th partial derivative. In case $k=0$, the 
$T_i$ reduce to the corresponding partial derivatives. In this paper, we shall
 assume throughout that $k\geq 0$ (i.e. all values of $k$ are non-negative), 
though several results of Section 3 may be extended to larger ranges of $k$.  
The most important basic properties of the $T_i$, proved in \cite{Du2}, are 
as follows: Let $\mathcal P=\b C\,[\b R^N]$ denote the algebra of polynomial 
functions on $\b R^N$ and $\mathcal P_n \>\, (n\in \b Z_+ = \{0,1,2\ldots\})$ 
the subspace of  homogeneous polynomials of degree $n$. Then \parskip=-1pt
\begin{enumerate}\itemsep=0pt
\item[(1.1)] The set $\{T_i\}$ generates a commutative algebra of 
differential-difference operators on $\mathcal P$.
\item[(1.2)] Each $T_i$ is homogeneous of degree $-1$ on $\mathcal P,$ that 
is, $T_i \,p \in \mathcal P_{n-1}$ for $p\in\mathcal P_n.$ 
\end{enumerate}
Of particular importance in this paper is   the generalized Laplacian 
associated with $G$ and $k$, which is defined as $\,\Delta_k := 
\sum_{i=1}^N T_i^2$. It is homogeneous of degree $-2$ on $\mathcal P$ and  
given  explicitely by
\[ \Delta_k f(x)\,=\, \Delta f(x) + 2\sum_{\alpha\in R_+} k(\alpha)
\Bigl[\frac{\langle\nabla f(x),\alpha\rangle}{\langle\alpha ,x\rangle} - 
\frac{f(x)-f(\sigma_\alpha x)}{\langle\alpha , x\rangle^2}\Bigr].\]
(Here $\Delta$ and $\nabla$ denote the usual Laplacian and gradient 
respectively). 

The operators $T_i$ were introduced and first studied by Dunkl in a series of 
papers ([Du1-4]) in connection with a generalization of  the classical theory 
of spherical harmonics: Here the uniform spherical surface measure on the 
(N-1)-dimensional unit sphere is modified by a weight function  which is 
invariant under the action of some finite reflection group $G$, namely
\[ w_k(x) = \prod_{\alpha\in R_+} |\langle\alpha,x\rangle|^{2k(\alpha)}\,, \]
where $k\geq 0$ is some fixed multiplicity function on the root system $R$ 
of $G$. Note that $w_k$ is homogeneous of degree $2\gamma$, with 
\[\gamma:=\sum_{\alpha\in R_+}k(\alpha).\]
 In this context, in \cite{Du3} the following bilinear form on $\mathcal P$ 
is introduced:
\[ [p,q]_k\,:=\,(p(T)\,q)(0)\>\>\text{ for }\, p,q\in \mathcal P\,.\]
Here $p(T)$ is the operator derived from $p(x)$ by replacing $x_i$ by $T_i$. 
Property (1.1) assures that $[.\,,.]_k$ is well-defined.  A useful collection 
of its properties can be found in \cite{DJO}. We just recall that $[.\,,.]_k$ 
is symmetric and positive-definite (in case $k\geq 0$), and that $[p,q]_k = 0$
 for $p\in \mathcal P_n,\, q\in \mathcal P_m$ with $n\not=m$. Moreover, 
$[.\,,.]_k$ is closely related to the scalar product on $L^2(\b R^N, 
w_k(x)e^{-|x|^2/2}dx)$: In fact, according to \cite{Du3}, \parskip=-1pt
\begin{enumerate}\itemsep=0pt
\item[(1.3)] $\displaystyle \qquad\qquad [p,q]_k\,=\, n_k\int_{\b R^N} 
e^{-\Delta_k/2} p(x)\, e^{-\Delta_k/2} q(x)\,w_k(x) e^{-|x|^2/2} dx \quad$ 
for all $p,q\in \mathcal P$, 
\end{enumerate}
with some normalization constant $n_k>0$. Given an orthonormal basis 
$\{\phi_\nu\,, \>\nu\in \b Z_+^N\}$ of $\mathcal P$ with respect to 
$[.\,,.]_k$, an easy rescaling of (1.3) shows that the polynomials    
\[ H_\nu(x) := 2^{|\nu|} e^{-\Delta_k/4}\varphi_{\nu}\]
are orthogonal  with respect to $\,w_k(x)e^{-|x|^2}dx$ on $\b R^N$. We call 
them the generalized Hermite polynomials on $\b R^N$ associated with $G,\, k$ 
and $\{\varphi_\nu\}$.  

\medskip

The first part of this paper is  devoted to the study of such Hermite 
polynomial systems and associated  Hermite functions.  They generalize their 
classical counterparts in a natural 
way: these are just obtained for $k=0$ and $\,\phi_\nu(x) = (\nu !)^{-1/2} 
x^\nu\,.$ In the one-dimensional case, associated with the reflection group 
$G=\b Z_2$ on $\b R$, our generalized Hermite polynomials coincide with those 
introduced in \cite{Ch} and studied in \cite{Ro}. Our setting also includes, 
for the symmetric group $G=S_N$,  the so-called non-symmetric generalized 
Hermite polynomials which were recently introduced by Baker and Forrester in 
\cite{BF2},  \cite{BF3}. These are non-symmetric analogues of Lassalle's  
(symmetric) generalized Hermite polynomials associated with the group $S_N$ 
(see \cite{La} and, for a further study, \cite{BF1}). Moreover, the 
``generalized Laguerre polynomials'' of \cite{BF2}, \cite{BF3} can be 
considered as a subsystem of Hermite polynomials associated with a reflection 
group of type $B_N$.

  After a short collection of notations and basic facts from Dunkl's theory in
 Section 2, the concept of generalized Hermite polynomials is introduced in 
Section 3, along with a discussion of the above mentioned special classes.  We
derive generalizations for many of the well-known  properties of the 
classical Hermite polynomials and Hermite functions: A Rodrigues formula, a 
generating relation and a Mehler formula for the Hermite polynomials, 
analoguesof the second order differential equations and a characterization of 
the generalized Hermite functions as eigenfunctions of the Dunkl transform. 
Parts of this section may be seen as a unifying treatment of results from 
\cite{BF2}, \cite{BF3} and \cite{Ro} for their particular cases.

In Section 4, which makes up the second major part of this paper, we turn to 
the Cauchy problem for the heat operator associated with the generalized 
Laplacian: Given an initial distribution $f\in C_b(\b R^N)$, there has to be 
found a function $u\in C^2(\b R^N\times (0,T))\cap C(\b R^N\times [0,T])$ 
satisfying
\begin{enumerate}\itemsep=0pt
\item[(1.4)] $\displaystyle \qquad\qquad
H_k u := \Delta_k u - \partial_t u\, = 0 \>\>$ on $\b R^N\times (0,\infty), 
\quad  \>u(.\,,0)=f.$
\end{enumerate}
For smooth and rapidly decreasing initial data $f$ an explicit solution is 
easy to obtain; it involves the generalized heat kernel
\[ \Gamma_k(x,y,t) = \frac{M_k}{t^{\gamma +N/2}}\,e^{-(|x|^2 + |y|^2)/4t}\,
K\Bigl(\frac{x}{\sqrt{2t}},\frac{y}{\sqrt{2t}}\Bigr),\quad x,y\in \b R^N,\> 
t>0.\]
Here $M_k$ is a positive constant and $K$ denotes the generalized exponential 
kernel associated with $G$ and $k$ as introduced in \cite{Du3}. In the theory 
of Dunkl operators and the Dunkl transform, it takes over the r\^{o}le of the 
usual exponential kernel $e^{\langle x,y\rangle}.$ Some of its properties are 
collected in Section 2. Without knowledge whether $K$ is nonnegative, a 
solution of (1.4) for arbitrary initial data seems to be difficult. However, 
one can prove a maximum principle for the generalized Laplacian $\Delta_k$, 
which is the key ingredient to assure that $\Delta_k$ leads to a positive 
one-parameter contraction semigroup on the Banach space $C_0(\b R^N), 
\|.\|_\infty)$.
Positivity of this semigroup enforces positivity of $K$  and allows to 
determine the explicit solution of (1.4) in the general case. We finish this 
section  with an extension of a well-known  maximum principle for the 
classical heat operator  to our situation.  This in particular implies a 
uniqueness result for solutions of the above Cauchy problem.  

\bigskip 

\noindent\emph{Acknowledgement.} It is a pleasure to thank Charles F. Dunkl 
and Michael Voit for some valuable comments and discussions.      

\section{Preliminaries}

The purpose of this section is  to  establish our basic notations and collect 
some further facts on Dunkl operators and the Dunkl transform which will be 
of importance later on. General references here are \cite{Du3}, \cite{Du4} 
and \cite{deJ}.

First of all we note the following product rule, which is confirmed by a 
short calculation:
\hfill\break
\noindent For each $f\in C^1(\b R^N)$ and each  $g\in C^1(\b R^N)$ which is 
invariant under the natural action of $G$, 
\begin{equation} T_i(fg)\,=\, (T_i f)g + f(T_i g) \quad\text{for } \>\, 
i=1,\ldots,N.
\end{equation}
We use the common multi-index notation; in particular, for  
$\nu=(\nu_1,\ldots,\nu_N)\in \b Z_+^N$ and $x=(x_1,\ldots,x_N)\in \b R^N$ we 
set $x^\nu := x_1^{\nu_1}\cdot\ldots\cdot x_N^{\nu_N}$, $\nu! := 
\nu_1!\cdot\ldots\cdot\nu_N!$  and $|\nu|:= \nu_1 +\ldots + \nu_N.$ If  
$f:\b R^N\to \b C$ is analytic with $f(x)=\sum_{\nu}a_\nu x^\nu\,,$  the 
operator $f(T)$ is defined by
\[f(T):= \sum_\nu a_\nu T^\nu\,=\, \sum_\nu a_\nu T_1^{\nu_1}\cdot\ldots\cdot 
T_N^{\nu_N}.\]
We restrict its action to $C^k(\b R^N)$ if $f$ is a polynomial of degree $k$ 
and to $\mathcal P$ otherwise. The following formula will be used frequently:

\begin{lemma} Let $p\in \mathcal P_n$. Then for $c\in \b C$ and 
$a\in \b C\setminus\{0\},$
\[ \bigl(e^{c\Delta_k}p\bigr)(ax) = a^n\bigl( e^{a^{-2}c\,\Delta_k}\bigr)
p(x)\quad \text{for all}\>\, x\in\b R^N.\] 
In particular, for $p\in\mathcal P_n$ we have
\begin{equation}\label{(1.30)}
\bigl( e^{-\Delta_k/2} p\bigr)(\sqrt 2 x)\,=\, \sqrt 2 ^{\,n} 
\bigl( e^{-\Delta_k/4} p\bigr)(x).
\end{equation}
\end{lemma}

\begin{proof} For $m\in \b Z_+$ with $2m\leq n,$ the polynomial 
$\Delta_k^mp$ is homogeneous of degree $n-2m.$ Hence
\[ \bigl(e^{c\Delta_k}p\bigr)(ax) = \sum_{m=0}^{\lfloor n/2\rfloor} 
\frac{c^m}{m!}(\Delta_k^m p)(ax)= \sum_{m=0}^{\lfloor n/2\rfloor}
\frac{c^m}{m!}a^{n-2m}(\Delta_k^mp)(x) = 
a^n\bigl( e^{a^{-2}c\,\Delta_k}p\bigr)(x).\]
\end{proof}

\cite{Du3} a

\medskip

A major tool in this paper is the  generalized exponential kernel $K(x,y)$ on 
$\b R^N \times \b R^N$, which generalizes the usual exponential function  
$e^{\langle x,y\rangle}$. It  was first introduced in \cite{Du3} by means of 
a certain intertwining operator. By a result of \cite{Op1} (see also 
\cite{deJ}), the function  $x\mapsto K(x,y)$  may be characterized  as the 
unique analytic solution of the system
$\, T_i f = y_i f\,\, (i=1,\ldots,N)$ on $\b R^N$ 
with  $f(0)=1.$  
Moreover, $K$ is symmetric in its arguments and  has a holomorphic extension 
to $\b C^N\times \b C^N$. Its power  series can be written as 
$\,K=\sum_{n=0}^\infty K_n$, where $K_n(x,y)= K_n(y,x)$ and $K_n$ is a 
homogeneous polynomial of degree $n$ in each of its variables. Note that 
$K_0=1$ and $K(z,0)= 1$ for all $z\in \b C^N$. 

For the reflection group $G=\b Z_2$ on $\b R,$ the multiplicity function $k$ 
is characterized by a single parameter $\mu \geq 0$, and  
the kernel $K$ is given explicitely by 
\[ K(z,w) = j_{\mu-1/2}(izw) - \frac{izw}{2\mu +1}\,j_{\mu+1/2}(izw),
\quad z,w\in \b C,\]
where for $\alpha\geq -1/2$, $j_\alpha$ denotes the normalized spherical 
Bessel function 
\[ j_\alpha(z) = 2^\alpha\Gamma(\alpha +1) \frac{J_\alpha(z)}{z^\alpha}\,=\, 
\Gamma(\alpha +1) \cdot\sum_{n=0}^\infty \frac{(-1)^n(z/2)^{2n}}{n! \,
\Gamma(n+\alpha +1)}\,.\]
For details and related material we refer to \cite{Du4}, \cite{Roe}, 
\cite{RV} and \cite{Ro}.

We list some further general properties of $K$ and the $K_n$ (all under the 
assumption  $k\geq 0$) from  \cite{Du3}, \cite{Du4} and \cite{deJ}:

\medskip
\noindent  For all $z,w\in \b C^N$ and $\lambda\in \b C$, 
\parskip=-3pt\nopagebreak
\begin{enumerate}\itemsep=0pt
\item[\rm{(2.3)}] $K(\lambda z,w) = K(z,\lambda w);$
\item[\rm{(2.4)}] $\displaystyle |K_n(z,w)| \leq  \frac{1}{n!} |z|^n|w|^n\,$  
and $\displaystyle \,|K(z,w)| \leq \,e^{|z||w|}\,.$
\end{enumerate}
For all $x,y\in \b R^N$ and $j=1,\ldots, N,$ \parskip=-1pt
\begin{enumerate}\itemsep=0pt
\item[\rm{(2.5)}] $ \>|K(ix,y)|\leq \sqrt{|G|}\,;$
\item[\rm{(2.6)}] $T_j^x K_n(x,y) = y_jK_{n-1}(x,y) \>\>\>\text{and} \>\>\>
T_j^x K(x,y) =  y_jK(x,y);$ \hfill\break
here the superscript  $x$ denotes that the operators act with respect to the 
$x$-variable.
\end{enumerate}\parskip=-2pt
In \cite{deJ}, exponential bounds for the usual partial derivatives of $K$ 
are given. They imply in particular that for each $\nu\in \b Z_+^N$ there 
exists a constant $d_\nu >0,$ such that \parskip=-3pt\nopagebreak
\begin{enumerate}\itemsep=0pt
\item[\rm{(2.7)}] $|\partial_x^\nu K(x,z)|\,\leq\, d_\nu\, |z|^{|\nu|}\,
e^{|x|\, |{\rm Re}\, z|}\>\> $ for all $x\in \b R^N, \> z \in \b C^N$.
\end{enumerate}\parskip=-1pt 
Let us finally recall a  useful reproducing kernel property of $K$ from 
\cite{Du4} (it is rescaled with respect to the original one, thus fitting 
better in our context of generalized Hermite polynomials): Define the 
probability measure $\mu_k$ on $\b R^N$ by
\[ d\mu_k(x) := c_k e^{-|x|^2}w_k(x)dx,\quad \text{ with }\>c_k = 
\Bigl(\int_{\b R^N} e^{-|x|^2} w_k(x) dx\Bigr)^{-1}.\] 
Moreover, for $z\in \b C^N$ set $\,l(z):= \sum_{i=1}^N z_i^2.$ Then for all 
$z,w\in \b C^N,$ 
\begin{enumerate}\itemsep=0pt
\item[\rm{(2.8)}] $\displaystyle  \int_{\b R^N} K(2z,x) K(2w,x)\,d\mu_k(x)\,
=\, e^{l(z)+l(w)} K(2z,w).$ 
\end{enumerate}

The generalized exponential function $K$ gives rise to an integral transform, 
called the Dunkl transform on $\b R^N,$ which was introduced in \cite{Du4} 
and has been thoroughly studied in \cite{deJ} for a large range of parameters 
$k$. The Dunkl transform associated with $G$ and $k\geq 0$ is defined by
\[ \mathcal D_k: L^1(\b R^N, w_k(x)dx)\to C(\b R^N);\>\>\>
 \mathcal D_k f(\xi):= \int_{\b R^N} f(x)\,K(-i\xi,x)\,w_k(x)dx \quad 
(\xi\in \b R^N).\]
In \cite{deJ}, many of the important properties of Fourier transforms on 
locally compact abelian groups are proved to hold true for $\mathcal D_k.$ In 
particular, 
$\mathcal D_k f\in C_0(\b R^N)$ for $f\in L^1(\b R^N, w_k(x)dx),$  
and there holds an $L^1$-inversion theorem, which we recall for later 
reference: If $f\in L^1(\b R^N, w_k(x)dx)$ with $\mathcal D_k f\in 
L^1(\b R^N, w_k(x)dx),$ then 
$\, f = 4^{-\gamma-N/2} c_k^2\,\mathcal E_k\mathcal D_k f \quad a.e.,\,$
where  $\mathcal E_k f(x) = \mathcal D_k f(-x).$ (Note that 
$\,\mathcal D_k(e^{-|x|^2/2})(0)\,=\, 2^{\gamma+N/2}c_k^{-1}$, which gives 
the connection of  our constant $c_k$  with that of de Jeu.) Moreover, 
the Schwartz space $\mathcal S(\b R^N)$ of rapidly decreasing functions on 
$\b R^N$ is invariant under $\mathcal D_k$, and $\mathcal D_k$ can be 
extended to a Plancherel transform on $L^2(\b R^N, w_k(x)dx)$.  For details 
see \cite{deJ}.

\section{Generalized Hermite polynomials and Hermite functions}

Let $\{\phi_\nu\,, \nu\in \b Z_+^N\}$ be an orthonormal basis of $\mathcal P$ 
with respect to the scalar product $[.\,,.]_k$ such that $\phi_\nu\in 
\mathcal P_{|\nu|}$ and the coefficients of the $\phi_\nu$ are real. As 
$\,\mathcal P = \bigoplus_{n\geq 0} \mathcal P_n$ and $\mathcal P_n \perp 
\mathcal P_m$ for $n\not=m$, the $\phi_\nu$ with $|\nu|=n$ can for example be 
constructed by Gram-Schmidt orthogonalization within $\mathcal P_n$ from an 
arbitrary ordered real-coefficient basis of   $\mathcal P_n$. If $k=0$, the 
Dunkl operator $T_i$ reduces to the usual partial derivative $\partial_i$, 
and the canonical choice of the basis $\{\phi_\nu\}$ is just $\,\phi_\nu(x):= 
(\nu!)^{-1/2} x^\nu.$

As in the classical case, we have the following  connection of the basis 
$\{\phi_\nu\}$ with the generalized exponential function $K$ and its 
homogeneous parts $K_n$:

\begin{lemma}
\begin{enumerate}\itemsep=0pt
\item[\rm{(i)}] $\displaystyle K_n(z,w) = \sum_{|\nu|=n} \phi_\nu(z)\,
\phi_\nu(w)\>$ for all $z,w\in \b C^N.$
\item[\rm{(ii)}] $\displaystyle K(x,y) = \sum_{\nu\in \b Z_+^N} \phi_\nu(x)\,
\phi_\nu(y)\>$ for all $x,y \in \b R^N$, \hfill\break
 where the convergence is absolute and locally uniform on 
$\b R^N\times \b R^N.$
\end{enumerate}
\end{lemma}

\begin{proof} (i) It suffices to consider the case $z,w\in \b R^N.$ So fix 
some $w\in \b R^N$. As a function of $z$, the polynomial $K_n(z,w)$ is 
homogeneous of degree $n$. Hence we have
\[ K_n(z,w) = \sum_{|\nu|=n} c_{\nu,\,w}\,\phi_\nu(z) \quad \text{with}\quad 
c_{\nu,w} = [K_n(.\,,w), \phi_\nu]_k\,.\]
Repeated application of formula (2.6) for $K_n$ gives
\[ c_{\nu,w} =  \phi_\nu(T^z)K_n(z,w)\, =\, \phi_\nu(w)\, K_0(z,w)\,=\, 
\phi_\nu(w).\]
 Thus part (i) is proved.  For (ii), first note that by  (2.4) we have   
$|K_n(x,x)| \leq \frac{1}{n!}|x|^{2n}$  and hence, as  the $\phi_\nu(x)$ are 
real,
$\,|\phi_\nu(x)|\leq \frac{1}{\sqrt{n!}}\,|x|^n\,$ for all $x\in \b R^N$ and 
all $\nu$ with $|\nu|=n.$
It follows that for  each $M>0$ 
the sum $\sum_{\b Z_+^N} |\phi_\nu(x)\phi_\nu(y)|\,$ is majorized on $\{(x,y):
 |x|,\,|y|\leq M\}$ by the convergent series $\,\sum_{n\geq 0} 
\binom{n+N-1}{n} M^{2n}/n!\,.$ This yields the assertion. 
\end{proof}

For homogeneous polynomials $p, q\in \mathcal P_n$, relation (1.3) can be 
rescaled (by use of formula \eqref{(1.30)}): 
\begin{equation}\label{(2.50)}
[p,q]_k\,=\, 2^n \int_{\b R^N} e^{-\Delta_k/4}p(x)\, e^{-\Delta_k/4}q(x)\, 
d\mu_k(x).
\end{equation}
This suggests to  define  a generalized multivariable Hermite polynomial 
system on $\b R^N$ as follows:

\begin{definition}
The generalized Hermite polynomials $\{H_\nu\,, \>\nu\in \b Z_+^N\}$  
associated with the basis $\{\phi_\nu\}$ on $\b R^N$ are given by 
\begin{equation}
H_\nu(x):= 2^{|\nu|}e^{-\Delta_k/4}\phi_\nu(x) = 
2^{|\nu|}\sum_{n=0}^{\lfloor|\nu|/2\rfloor} \frac{(-1)^n}{4^n n!}\,
\Delta_k^n \phi_\nu(x).
\end{equation}
Moreover, we define the generalized Hermite functions on $\b R^N$ by
\begin{equation}
h_\nu(x):= e^{-|x|^2/2}H_\nu(x), \quad \nu\in \b Z_+^N.
\end{equation}
\end{definition}

\medskip
Note that $H_\nu$ is a polynomial of degree $|\nu|$ satisfying $H_\nu(-x) = 
(-1)^{|\nu|}H_\nu(x)$ for all $x\in \b R^N$.  A standard argument shows  that 
$\mathcal P$ is dense in $L^2(\b R^N, d\mu_k)$. Thus by virtue  of 
\eqref{(2.50)}, the $\{2^{-|\nu|/2}H_\nu, \>\nu\in \b Z_+^N\}$ form an 
orthonormal basis of $L^2(\b R^N, d\mu_k).$ Let us give two immediate examples:

\medskip
\begin{examples}
\begin{enumerate}\itemsep=1pt
\item[\rm{(1)}] In the classical case $k=0$ and  $\,\phi_\nu(x):= 
(\nu!)^{-1/2} x^\nu,$  we obtain
\[ H_\nu(x) = \frac{2^{|\nu|}}{\sqrt{\nu!}}\prod_{i=1}^N e^{-\partial_i^2/4}
(x_i^{\nu_i})\,=\,\frac{1}{\sqrt{\nu!}}\prod_{i=1}^N \widehat H_{\nu_i}(x_i),\]
where  the $\widehat H_n,\>n\in \b Z_+$ denote the classical Hermite 
polynomials on $\b R$ defined by 
\[ e^{-x^2} \widehat H_n(x) = (-1)^n \,\frac{d^n}{dx^n}\bigl(e^{-x^2}\bigr).\]
\item[\rm{(2)}] For the reflection group $G=\b Z_2$ on $\b R$ and multiplicity
  parameter $\mu \geq 0$,  
 the polynomial basis $\{\phi_n\}$ on $\b R$ with respect to $[.\,,.]_\mu$ is 
determined uniquely (up to sign-changes) by suitable normalization of the 
monomials $\{x^n, n\in\b Z_+\}$. One obtains $\,H_n(x) = d_n H_n^\mu(x),$ 
where $d_n\in \b R\setminus\{0\}$ are constants and the $H_n^\mu \,,\> n\in 
\b Z_+$ are the generalized Hermite polynomials on $\b R$ as introduced e.g. 
in \cite{Ch} and studied in \cite{Ro} (in some different normalization). They 
are orthogonal with respect to $|x|^{2\mu} e^{-|x|^2}$ and can be written  as
\[ \begin{cases}
     H_{2k}^\mu(x) = (-1)^k 2^{2k} k! \,L_k^{\mu-1/2}(x^2),\\
     H_{2k+1}^\mu(x) = (-1)^k 2^{2k+1} k! \,xL_k^{\mu+1/2}(x^2);
   \end{cases}\]
here  the $L_n^\alpha$ are the Laguerre polynomials of index 
$\alpha\geq -1/2$,  given by
\[ L_n^\alpha(x) = \frac{1}{n!} x^{-\alpha} e^x\,\frac{d^n}{dx^n}
\Bigl( x^{n+\alpha} e^{-x}\Bigr).\]
\end{enumerate}
\end{examples} 

\parskip=1pt

Before discussing further examples, we are going to establish generalizations 
of the classical second order differential equations for Hermite polynomials 
and Hermite functions. For their proof we shall employ the 
$sl(2)$-commutation relations of the operators
\[ E:= \frac{1}{2}|x|^2,\>\> F:= -\frac{1}{2}\Delta_k \quad\text{and}\>\> 
H:= \sum_{i=1}^N x_i\partial_i +\left(\gamma +N/2\right)\]
on $\mathcal P$, which can be found e.g. in \cite{He}; they are 
\begin{equation}\label{(2.40)}
\big[H,E\big] = 2E,\>\> \big[H,F\big] = -2F,\>\>\big[E,F\big] = H.
\end{equation}
(As usual, $\,\big[A,B\big] = AB-BA$ for operators $A,B$ on $\mathcal P$.) 
The first two relations are immediate consequences of the fact that the 
Euler operator  $\,\rho:=\sum_{i=1}^N x_i\partial_i\,$ satisfies  
$\,\rho(p)=np\,$ for each homogeneous $p\in \mathcal P_n$. We have the 
following general result:

\begin{theorem}
\begin{enumerate}
\item[\rm{(1)}] For $n\in \b Z_+$ set $\,V_n:= \{e^{-\Delta_k/4} p\,: 
p\in \mathcal P_n\}.$ Then $\,\mathcal P = \bigoplus_{n\in \b Z_+} V_n$, 
and $V_n$ is the eigenspace of the operator $\,\Delta_k -2\rho\,$ on 
$\mathcal P$ corresponding to the eigenvalue $-2n$.
\item[\rm{(2)}] For $q\in V_n$, the function $\,f(x):= e^{-|x|^2/2} q(x)$ 
satisfies
\[ \bigl(\Delta_k -|x|^2\bigr) f \,=\, -(2n + 2\gamma+N)f\,.\]
\end{enumerate}
\end{theorem} 

\begin{proof}
(1) It is clear that $\,\mathcal P = \bigoplus V_n\,$. By induction from 
\eqref{(2.40)}, or simply by the facts that $\rho(p)=np$ for 
$p\in\mathcal P_n$ and $\Delta_k$ is homogeneous of degree $-2$,  
we obtain the commuting relations 
\[\big[2\rho,\Delta_k^n\big] = -4n\Delta_k^n\quad\text{ for all}\>\>
n\in \b Z_+\,,\quad\text{hence} \quad  
\big[ 2\rho, e^{-\Delta_k/4}\big] = \Delta_k e^{-\Delta_k/4}.\] 
For arbitrary $q\in \mathcal P$ and $p:= e^{\Delta_k/4}q\,$ it now follows 
that 
\[2\rho(q) \,=\, (2\rho e^{-\Delta_k/4})(p)\, =\, 2e^{-\Delta_k/4}\rho(p) + 
\Delta_k e^{-\Delta_k/4}p\,=\, 2 e^{-\Delta_k/4}\rho(p) +\Delta_k q.\]
Hence for $a\in \b C$, there are equivalent:
\[ (\Delta_k -2\rho)(q) = -2aq \,\Longleftrightarrow \, \rho(p) = ap\,
\Longleftrightarrow\, a=n\in \b Z_+ \>\>\text{and}\>\> p\in \mathcal P_n.\]
This yields the assertion.

(2) We first verify by induction that
\begin{equation}\label{(2.41)}
 \big[\Delta_k, E^n\big] = 2nE^{n-1}H + 2n(n-1) E^{n-1}\quad\text{for all}\>\>
 n\in \b N.
\end{equation}
In case $n=1$ this is clear from \eqref{(2.40)}; if $n\geq 1$, then
\[\big[\Delta_k, E^{n+1}\big] = \big[ \Delta_k, E^n\big]\,E + 
E^n\big[\Delta_k,E\big] = 2nE^{n-1}HE + 2n(n-1)E^n + 2E^nH,\]
by our induction hypothesis. Using the identity $\,HE = EH +2E$, we readily 
see that \eqref{(2.41)} holds for $n+1$ as well. 

From \eqref{(2.41)} it is now easily deduced that $\big[\Delta_k, e^{-E}\,
\big] = -2e^{-E}H + 2Ee^{-E}.$ We thus obtain
\begin{align}
(\Delta_k - |x|^2)f\, = \,\Delta_k\bigl(e^{-E} q\bigr) -2E e^{-E} q\, = 
\,e^{-E}\Delta_k q -2 e^{-E}\,(\rho +\gamma +N/2) q\,.
\end{align}
The stated relation is now a consequence of (1).

\end{proof}

\begin{corollary} 
\begin{enumerate}\itemsep=0pt
\item[\rm{(i)}] The generalized Hermite polynomials  satisfy the following 
differential- difference  equation:
\[\Bigl(\Delta_k - 2\sum_{i=1}^N x_i\partial_i\Bigr) H_\nu\,=\, 
-2|\nu|H_\nu\,, \quad\> \nu\in \b Z_+^N.\]
\item[\rm{(ii)}] The generalized Hermite functions $\{h_\nu\,,\nu\in 
\b Z_+^N\}$ form a complete set of eigenfunctions for the operator 
$\,\Delta_k -|x|^2\,$ on $L^2(\b R^N, w_k(x)dx)$ with 
\[ \bigl(\Delta_k - |x|^2\bigr)\, h_\nu\,=\, -(2|\nu| +2\gamma +N)\,h_\nu\,.\]
 \end{enumerate}
\end{corollary}

Note also that as a consequence of the above theorem, the operator 
$\,\Delta_k-2\rho\,$ has for each $p\in \mathcal P_n$ a unique polynomial 
eigenfunction $q$ of the form $\,q=p+r$, where the degree of $r$ is less 
than $n$; it is  given by $q = e^{-\Delta_k/4} p.$

\begin{examples}
\begin{enumerate}
\item[\rm{(3)}] {\bf The $S_N$-case.}\enskip 
For the symmetric group  $G=S_N$ (acting on $\b R^N$ by permuting the 
coordinates), the mulitplicity function is characterized by a single 
parameter which is often denoted by $1/\alpha>0$,  and the corresponding 
weight function is given by $\,w_S(x)= \prod_{i<j} |x_i-x_j|^{2/\alpha}.$ 
The associated Dunkl operators are 
\[ T_i^S \,=\,\partial_i + \frac{1}{\alpha} \sum_{j\not= i} 
\frac{1-s_{ij}}{x_i-x_j} \quad (i=1,\ldots,N),\]
where $s_{ij}$ denotes the operator transposing $x_i$ and $x_j$. 
The operator $\Delta_S -2\rho$ is a Schr\"odinger operator of 
Calogero-Sutherland type, involving an external harmonic potential with 
exchange terms, see \cite{BF2} and \cite{BF3}. It is given explicitely by
\begin{equation}\label{(1.45)}
 \Delta_S -2\rho\,=\,\Delta -2\sum_{i=1}^N x_i\partial_i + \frac{2}{\alpha} 
\sum_{i<j} \frac{1}{x_i-x_j} \Big[(\partial_i -\partial_j) - 
\frac{1-s_{ij}}{x_i-x_j}\Big]\,.
\end{equation}
In \cite{BF2}, Baker and Forrester study ``non-symmetric generalized Hermite 
polynomials''  $E_\nu^{(H)},$ which they define as the unique eigenfunctions 
of \eqref{(1.45)} of the form 
\[ E_\nu^{(H)}\,=\, E_\nu + \sum_{|\mu|<|\nu|}c_{\mu,\,\nu} E_\mu\,,\]
where the $E_\nu,\> \nu\in \b Z_+^N\,$ are the non-symmetric Jack polynomials 
(associated with $S_N$ and $\alpha$) as defined e.g. in \cite{Op} (see also 
\cite{KS}). Thus $\, E^{(H)}_\nu = e^{-\Delta_S/4} E_\nu$ (just by 
Lemma 3.4), and indeed, up to some normalization factors, the $E^{(H)}_\nu$ 
make up a system of generalized Hermite polynomials for $S_N$ in our sense. 
This follows from the fact that the non-symmetric Jack polynomials $E_\nu$, 
being homogeneous of degree $|\nu|$ and forming a vector space basis of 
$\mathcal P$, are also orthogonal with respect to the inner product 
$[.\,,.]_S$. This  was proved in \cite{BF3} via  orthogonality of the 
$E_\nu^{(H)}$. A short direct proof can be given as follows: According to 
\cite{Op}, Prop. 2.10, the $E_\nu$ are simultaneous eigenfunctions of the 
Cherednik operators $\xi_i$ for $S_N$, which were introduced in \cite{Che} 
and  can be written as 
\begin{equation}\label{(2.51)}
 \xi_i = \alpha x_i T_i^S + 1-N + \sum_{j>i} s_{ij} \qquad (i=1,\ldots,N).
\end{equation}
In fact, the $E_\nu$ satisfy $\, \xi_i E_\nu = \overline\nu_i E_\nu\,, $ 
where the eigenvalues $\overline\nu = (\overline\nu_1,\ldots,\overline\nu_N)$ 
are given explicitely in \cite{Op}. They are distinct, i.e. if $\nu\not=\mu$, 
then $\overline\nu\not=\overline\mu$. 
On the other hand, \eqref{(2.51)}  shows immediately that the operators 
$\xi_i$ are symmetric with respect to $[.\,,.]_S$. (For all $p,q\in 
\mathcal P$ and $g\in S_N$, we have $\,[g(p), g(q)]_S = [p,q]_S$ for 
$g(p)(x)= p(g^{-1}(x))$,  as noted in \cite{DJO}.) Together, this proves 
that the $E_\nu$ are orthogonal with respect to $[.\,,.]_S\,$. Hence a 
possible choice for the basis $\{\varphi_\nu\}$ is to set $\phi_\nu = d_\nu 
E_\nu$, with some normalization constants $d_\nu>0$. 

We finally remark that in this case the locally uniform convergence of the 
series in Lemma 3.1(ii)  extends to $\b C^N\times \b C^N$, see also  
\cite{BF3}, Prop. 3.10. This is because     the coefficients of the 
non-symmetric Jack-polynomials $E_\nu$ in their monomial expansions  are 
known to be nonnegative (\cite{KS}, Theorem 4.11), hence $|E_\nu(z)|\leq 
E_\nu(|z|)$ for all $z\in \b C^N$.  

\item[\rm{(4)}] {\bf A remark on the $B_N$-case.}\enskip Suppose that $G$ is 
the Weyl group of type $B_N$, generated by sign-changes and permutations. 
Here the multiplicity function is characterized by two parameters $k_0,\, 
k_1\geq 0.$
The weight function is 
\[w_B(x) = \prod_{i=1}^N |x_i|^{2k_1} \prod_{i<j} |x_i^2-x_j^2|^{2k_0}.\]
  Let $T_i^B$ and $\Delta_B$ denote the associated Dunkl operators and 
Laplacian. We consider the space 
\[ W:= \{ f\in C^1(\b R^N): \,f(x) = F(x^2) \quad\text{for some} \> 
F\in C^1(\b R^N)\}\]
of ``completely even'' functions;
here $\,x^2 = (x_1^2,\ldots, x_N^2).$ It is easily checked that for 
completely even $f$, $\Delta_B f$ is also completely even. The restriction 
of $\Delta_B$ to $W$ is given by
\begin{align} 
\Delta_B\vert_W \,=\,& \Delta + 2k_1\sum_{i=1}^N \frac{1}{x_i}\,
\partial_i\, +\, 2k_0 \sum_{i<j} \Bigl(\frac{1}{x_i -x_j}
\bigl(\partial_i -\partial_j\bigr) + \frac{1}{x_i + x_j}
\bigl(\partial_i + \partial_j\bigr)\Bigr) \notag\\
-\,& 2k_0 \sum_{i<j} \Bigl( \frac{1}{(x_i-x_j)^2} +\frac{1}{(x_i+x_j)^2}\Bigr)
 \bigl(1-s_{ij}\bigr). \notag\end{align}
Again, the operator $(\Delta_B - 2\rho)\vert_W\,$ is of Calogero-Sutherland 
type. Its completely even polynomial eigenfunctions are discussed in 
\cite{BF2} and \cite{BF3}  separately from the Hermite-case; they are 
called ``non-symmetric Laguerre polynomials'' and denoted  by 
$E_\nu^{(L)}(x^2).$ It is easy to see that they make up the completely even 
subsystem of a suitably choosen  generalized Hermite-system $\{H_\nu\}$ for 
$B_N$ (and parameters $k_0,\,k_1$, where we assume $k_0>0$):

 To this end, let again
 $E_\nu$ denote the $S_N$-type non-symmetric Jack polynomials, corresponding 
to $\alpha =1/k_0$. For $\nu\in \b Z_+^N$ set 
$\widehat E_\nu(x):= E_\nu(x^2).$  These modified Jack polynomials form a 
basis of $\mathcal P\cap W$. The non-symmetric Laguerre polynomials of Baker 
and Forrester can be written as 
\[E_\nu^{(L)}(x^2)\,=\,e^{-\Delta_B/4}\widehat E_\nu (x)\,.\]  
(Note that the polynomials on the right side are in fact completely even and 
eigenfunctions of $\Delta_B -2\rho$.) Involving again the $S_N$-type 
Cherednik operators from (3), it is easily checked that the 
$\,\widehat E_\nu$ are orthogonal with respect to the scalar product 
$[.\,,.]_B$: The $\xi_i$ induce operators $\widehat\xi_i \>(i=1,\ldots,N)$ 
on $W$ by
\[\, \widehat \xi_i f(x) := (\xi_iF)(x^2) \quad\text{if }\> f(x) = F(x^2),\] 
c.f. \cite{BF3}. Thus $\,\widehat \xi_i \widehat E_\nu = \overline\nu_i 
\widehat E_\nu$, and a short calculation gives
\[ \widehat\xi_i f(x)\,=\, \alpha x_i^2 (T_i^S F)(x^2) + 1-N + \sum_{j>i} 
s_{ij}\,=\, \frac{\alpha}{2} x_i T_i^B f(x) + 1-N + \sum_{j>i} s_{ij}.\]
This proves that the $\widehat \xi_i$ are symmetric with respect to 
$[.\,,.]_B$ and yields our assertion by the same argument as in the previous 
example. We therefore obtain an orthonormal  basis $\{\varphi_\nu\}$ of 
$\mathcal P$ with respect to $[.\,,.]_B$ by setting $\,\varphi_\nu := 
d_\nu\widehat E_\eta$ for $\nu=(2\eta_1,\ldots,2\eta_N)$ and completing 
the set $\{\varphi_\nu\,, \> \nu\in (2\b Z_+)^N\}$ by a Gram-Schmidt procedure.

\end{enumerate}
\end{examples}

Many properties of the classical Hermite polynomials and Hermite functions on 
$\b R^N$ have natural extensions to  our  general setting. We start with a 
Rodrigues-formula:

\begin{theorem} For all $\nu\in \b Z_+^N$ and $x\in \b R^N$ we have
\begin{equation}
H_\nu(x) = (-1)^{|\nu|}e^{|x|^2}\phi_\nu(T)\,e^{-|x|^2}.
\end{equation}
\end{theorem}

\begin{proof}
First note that if $p$ is a polynomial of degree $n\geq 0$, then 
\[p(T)\, e^{-|x|^2} = q(x)\,e^{-|x|^2}\]
with a polynomial $q$ of the same degree. This follows easily from induction 
by the degree of $p$, together with the product rule (2.1). In particular, 
the function 
\[ Q_\nu(x):= (-1)^{|\nu|} e^{|x|^2}\phi_\nu(T)\, e^{-|x|^2} = 
e^{|x|^2}\phi_\nu(-T)\,e^{-|x|^2}\]
is a polynomial of degree $|\nu|.$ In order to prove that $Q_\nu = H_\nu$, 
it therefore suffices to show that for each $\eta\in \b Z_+^N$ with 
$|\eta|\leq|\nu|$, 
\begin{equation}\label{(2.1)}
2^{-|\eta|}\,\int_{\b R^N} Q_\nu(x) H_\eta(x)\,d\mu_k(x)\,=\,
\delta_{\nu,\eta}\,,
\end{equation}
where $\,\delta_{\nu,\eta}\,$ denotes the Kronecker delta. Using the 
antisymmetry of the $\,T_i\,$ with respect to $\,L^2(\b R^N, w_k(x)dx)$ 
(Lemma 2.9 of \cite{Du4}) as well as the commutativity of $\{T_i\},$ we can 
write
\begin{align}
 2^{-|\eta|}\int_{\b R^N} &Q_\nu(x) H_\eta(x)\,d\mu_k(x) = c_k\int_{\b R^N} 
\phi_\nu(-T)\bigl(e^{-|x|^2}\bigr) e^{-\Delta_k/4}\phi_\eta(x)\,w_k(x)dx 
\notag  \\ \notag
& = c_k\int_{\b R^N} e^{-|x|^2}\bigl(\phi_\nu(T) 
e^{-\Delta_k/4}\phi_\eta\bigr)(x)\,w_k(x)dx = \int_{\b R^N} 
\bigl(e^{-\Delta_k/4}\phi_\nu(T)\,\phi_\eta\bigr)(x)\,d\mu_k(x).  \notag
\end{align}
But as $|\eta|\leq |\nu|$,  we have $\,\phi_\nu(T)\,\phi_\eta = 
[\phi_\nu, \phi_\eta]_k = \delta_{\nu, \eta}\,$   from which \eqref{(2.1)} 
follows.
\end{proof}

There is also a generating function for the generalized Hermite polynomials:

\begin{proposition}
For $n\in \b Z_+$ and $z,w\in \b C^N$ put 
$\, L_n(z,w):= \sum_{|\nu|=n} H_\nu(z)\, \phi_\nu(w).\,$
Then
\[\sum_{n=0}^\infty L_n(z,w)\,=\, e^{-l(w)} K(2z,w),\]
the convergence of the series being normal on $\b C^N\times \b C^N$.
\end{proposition}

\begin{proof}
Suppose first that $z,w\in\b R^N$. 
By definition of the $H_\nu$ and in view of formula (2.6) for $K_n$ we may 
write
\begin{align}
 L_n(z,w) =\, &2^n e^{-\Delta_k^z/4} K_n(z,w) = 
2^n\sum_{m=0}^{\lfloor n/2\rfloor} \frac{(-1)^m}{4^mm!}l(w)^{m} 
K_{n-2m}(z,w)\notag \\
 = \,&\sum_{m=0}^{\lfloor n/2\rfloor} \frac{(-1)^m}{m!} l(w)^m 
K_{n-2m}(2z,w)\notag
\end{align}
for all $n\in \b Z_+$. By analytic continuation, this holds for all 
$z,w\in \b C^N$ as well. Using estimation (2.4), one obtains 
\[ S_n(z,w):= \sum_{m=0}^{\lfloor n/2\rfloor}\frac{1}{m!}\,
|l(w)|^{m}|K_{n-2m}(2z,w)|\,  \leq \sum_{m=0}^{\lfloor n/2\rfloor}
\frac{1}{m!}\,|w|^{2m}\cdot\frac{|2z|^{n-2m}|w|^{n-2m}}{(n-2m)!}.\]
If $n$ is even, set $k:=n/2$ and estimate further as follows:
\[ S_n(z,w)\leq \frac{|w|^{2k}}{ k!}\sum_{m=0}^k \binom{k}{m}(2|z|^2)^{k-m}
\,=\, \frac{1}{k!}\bigl(|w|^2(1+2|z|^2)\bigr)^k.\]
A similar estimation holds if  $n$ is odd. 
This entails  the normal convergence of the series \hfill\break
$\sum_{n=0}^\infty L_n(z,w)$ on $\b C^N\times \b C^N$, and also that
\begin{align}
 \sum_{n=0}^\infty L_n(z,w)\,&=\sum_{n=0}^\infty\sum_{m=0}^\infty 
\frac{(-1)^m}{m!}\,l(w)^{m}K_{n-2m}(2z,w) \quad\text{(with }\> \,
K_j:=0\>\>\text{for }\>\, j<0)\notag \\
&= \sum_{m=0}^\infty \frac{(-1)^m}{ m!}\,l(w)^{m}\sum_{n=0}^\infty 
K_{n-2m}(2z,w)\,=\, e^{-l(w)} K(2z,w) \notag 
\end{align}
for all $z,w\in \b C^N$.
\end{proof}

Applying Lemma 2.1  to $p=\phi_\nu$ with $c=-1/4$ and $a=1/\lambda$, we 
obtain the following formula for the generalized Hermite polynomials:

\begin{lemma}
For $\lambda\in \b C\setminus \{0\},\> \nu\in \b Z_+^N$ and $x\in \b R^N$, 
\[
\Bigl(\frac{\lambda}{2}\Bigr)^{|\nu|} H_\nu\Bigl(\frac{x}{\lambda}\Bigr)\,=\,
 \bigl( e^{-\lambda^2\Delta_k/4}\phi_\nu\bigr)(x).
\]
\end{lemma}

\begin{proposition} The generalized Hermite functions $\{h_\nu,\, \nu\in 
\b Z_+^N\}$ are a basis of  eigenfunctions of the Dunkl transform 
$\mathcal D_k$ on $L^2(\b R^N, w_k(x)dx)$, satisfying
\[ \mathcal D_k(h_\nu)\,=\, 2^{\,\gamma+N/2} c_k^{-1}\,(-i)^{|\nu|}\,h_\nu\,.\]
\end{proposition}

\begin{proof} We use Prop. 2.1 form \cite{Du4}, which says that for all 
$p\in \mathcal P$ and $z\in \b C^N$, 
\begin{equation}\label{(2.60)}
\frac{c_k}{2^{\gamma+N/2}}\int_{\b R^N} e^{-\Delta_k/2} p(x)\,K(x,z)\,w_k(x) 
e^{-|x|^2/2} dx\,=\, e^{l(z)/2} p(z).
\end{equation}
Here again, $\,l(z) = \sum_{i=1}^N z_i^2.$ Let $p_\nu(x):= e^{\Delta_k/2} 
H_\nu(x).$ In view of \eqref{(2.60)} we can write
\[ \mathcal D_k(h_\nu) (\xi) = \int_{\b R^N}  H_\nu(x) K(-i\xi,x)\, w_k(x) 
e^{-|x|^2/2}dx\, =\, 2^{\gamma+N/2} c_k^{-1}\, e^{-|\xi|^2/2}p_\nu(-i\xi)\]
for all $\xi\in \b R^N$. By definition of $H_\nu$ we have  $\, p_\nu(x) = 
2^{|\nu|} e^{\Delta_k/4} \varphi_\nu(x).$ So we arrive at
\[\mathcal D_k(h_\nu)(\xi) = 2^{\gamma+N/2} c_k^{-1}\, e^{-|\xi|^2/2}\, 
2^{|\nu|}\bigl(e^{\Delta_k/4}\phi_\nu\bigr)(-i\xi).\]
Application of Lemma 3.9  with $\lambda=-i$ now yields that $\, 
\bigl(e^{\Delta_k/4} \varphi_\nu\bigr)(-i\xi) = (-i/2)^{|\nu|} H_\nu(\xi)$, 
hence
\[ \mathcal D_k(h_\nu)(\xi)\, =\, 2^{\gamma+N/2} c_k^{-1} (-i)^{|\nu|} 
h_\nu(\xi).\]
\end{proof}

We finish this section with a Mehler-type formula for the generalized 
Hermite polynomials. For this, we need the following integral representation:

\begin{lemma}
For all $x,y\in \b R^N$ and $\nu\in \b Z_+^N$ we have
\[ e^{-|x|^2} H_\nu(x) \,=\,2^{|\nu|}\int_{\b R^N} K(x,-2iy)\,\phi_\nu(iy)\,
d\mu_k(y).\]
\end{lemma}

\begin{proof} A short calculation, using again relation  (2.2), shows that 
for homogeneous polynomials $p$,  formula \eqref{(2.60)} may be rewritten as
\begin{equation}\label{(2.65)}
\int_{\b R^N} e^{-\Delta_k/4} p(x) \,K(x,2z)\,d\mu_k(x) \,=\, e^{l(z)} p(z) 
\quad (z\in \b C^N).
\end{equation}
By linearity, this holds for all $p\in\mathcal P.$ Lemma 3.9 with $\lambda=i$ 
further shows that 
\[e^{\Delta_k/4}\phi_\nu(x) = \bigl(\frac{i}{2}\bigr)^{|\nu|} H_\nu(-ix).\]
 As $\varphi_\nu$ is homogeneous of degree $|\nu|,$ we thus can write  
$\phi_\nu(2iy) = (-i)^{|\nu|}\, e^{-\Delta_k/4}H_\nu^*(y)$ with $H_\nu^*(y) 
= H_\nu(iy).$ From \eqref{(2.65)}  it now follows that
\[\int_{\b R^N} K(x,-2iy)\, \phi_\nu(2iy)\, d\mu_k(y) \,= \,e^{-|x|^2} 
H_\nu^*(-ix),\]
which yields the assertion.
\end{proof}

\begin{theorem} {\rm (Mehler-formula for the $H_\nu$)} For $r\in \b C$ with 
$|r|<1$ and all $x,y \in \b R^N,$ 
\[ \sum_{\nu\in \b Z_+^N} \frac{H_\nu(x) H_\nu(y)}{2^{|\nu|}}\, r^{|\nu|} \,=
\, \frac{1}{(1-r^2)^{\,\gamma+N/2}}\,
{\rm exp}\left\{-\frac{r^2(|x|^2 + |y|^2)}{1-r^2}\right\} 
K\!\left(\frac{2rx}{1-r^2},\,y\right).\]
\end{theorem}

\begin{proof}
Consider the integral 
\[ M(x,y,r):= c_k^2 \int_{\b R^N\times\b R^N} K(-2rz,v) K(-2iz,x) K(-2iv,y)\, 
w_k(z) w_k(v)\, e^{-(|z|^2+ |v|^2)} d\,(z,v).\]
The bounds (2.4) and (2.5) on $K$ assure that it converges for all 
$r\in \b C$ with $|r|<1$ and all $x,y~\in~\b R^N.$ Now write
$\, K(-2rz,v)= \sum_{n=0}^\infty (2r)^n K_n(iz,iv)$ in the integral above. As
\[\sum_{n=0}^\infty |2r|^n |K_n(iz,iv)|\leq e^{2|r||z||v|}\,,\] 
the dominated convergence theorem  yields that
\begin{align}
 M(x,y,r) =& \,\sum_{n=0}^\infty (2r)^n \int_{\b R^N}\int_{\b R^N} 
K_n(iz,iv)\, K(-2iz,x)\, K(-2iv,y)\,d\mu_k(z)\, d\mu_k(v)  \notag\\
=& \,\sum_{n=0}^\infty (2r)^n \sum_{|\nu|=n} \Bigl(\int_{\b R^N} K(-2iz,x)\, 
\phi_\nu(iz)\, d\mu_k(z)\Bigr)\Bigl( \int_{\b R^N} K(-2iv,y)\,\phi_\nu(iv)\,
d\mu_k(v)\Bigr).\notag  \end{align}
From the above lemma we thus obtain
\begin{equation}\label{(2.20)}
 M(x,y,r) =\, e^{-(|x|^2 +|y|^2)} \sum_{\nu\in \b Z_+^N} r^{|\nu|} 
\frac{H_\nu(x)H_\nu(y)}{2^{|\nu|}},
\end{equation}
and this series, as a power series in $r,$ converges absolutely for all 
$x,y\in \b R^N.$ On the other hand, iterated integration and repeated 
application of formula (2.3) and the reproducing formula (2.8)  show that 
for real $r$ with $|r|<1$ we have 
\begin{align}
 M(x,y,r) =&\,c_k\int_{\b R^N}\Bigl(\int_{\b R^N} K(-2rz,v)\,K(-2iy,v) 
\,d\mu_k(v)\Bigr) K(-2iz,x) e^{-|z|^2} w_k(z)dz \notag \\
=&\, c_k\, e^{-|y|^2} \int_{\b R^N} e^{(r^2-1)|z|^2} K(2iry,z)\, K(-2ix, z)\, 
w_k(z)\,dz \notag\\
 =&\, c_k (1-r^2)^{-(\gamma +N/2)} e^{-|y|^2} \int_{\b R^N} e^{-|u|^2} 
K\Bigl( u\,, \frac{2iry}{\sqrt{1-r^2}}\Bigr) K\Bigl(u\,,
\frac{-2ix}{\sqrt{1-r^2}}\Bigr) w_k(u)\, du \notag\\
=&\, (1-r^2)^{-(\gamma +N/2)}\, 
{\rm exp}\!\left\{-\frac{|x|^2 + |y|^2}{1-r^2}\right\} 
K\Bigl(\frac{2rx}{1-r^2},\,y\Bigr).\notag
\end{align}
By analytic continuation, this holds for  $\{r\in \b C: |r|<1\}$ as well. 
Together with \eqref{(2.20)}, this finishes the proof.
\end{proof}

\section{The heat equation for Dunkl operators}

As before, let $\Delta_k$ denote the generalized Laplacian associated with 
some finite reflection group $G$ on $\b R^N$ and a multiplicity function 
$k\geq 0$ on its root system $R$. Recall that its action on $C^2(\b R^N)$ is 
given by 
\[\Delta_k f \,=\, \Delta f + 2\sum_{\alpha\in R_+} k(\alpha)\,\delta_\alpha 
f\,,\]
where
\[ \delta_\alpha f(x)\,=\, 
\frac{\langle\nabla f(x),\alpha\rangle}{\langle\alpha ,x\rangle} - \,
\frac{f(x)-f(\sigma_\alpha x)}{\langle\alpha , x\rangle^2}\,.\]
Its action may as well be restricted to $C^2(\Omega)$, where 
$\Omega\subset \b R^N$ is open and invariant under the group operation of 
$G$. We call a function $f\in C^2(\Omega)$  $k$-subharmonic on $\Omega,$  if 
$\Delta_k f\geq 0$ on $\Omega$. 

The generalized Laplacian satisfies the following maximum principle, which 
will be important later on:

\begin{lemma}
Let $\Omega\subseteq\b R^N$ be open and $G$-invariant. 
If a real-valued function $f\in C^2(\Omega)$ attains an absolute maximum at  
$x_0\in\Omega$, i.e. $f(x_0)= \sup_{x\in\Omega} f(x)$, then
\[\Delta_k f(x_0)\,\leq\, 0\,.\]
\end{lemma}

\begin{proof}
Let $D^2f(x)$ denote the Hessian of $u$ in $x\in \Omega$. 
The given situation enforces that $\nabla f(x_0) =0$ and $D^2f(x_0)$ is 
negative semi-definite; in particular, $\Delta f(x_0) \leq 0.$ Moreover, 
$f(x_0) \geq f(\sigma_\alpha x_0)$ for all $\alpha\in R$, so the statement is 
obvious in the case that $\langle\alpha ,x_0\rangle \not=0$ for all 
$\alpha\in R$. If $\langle\alpha ,x_0\rangle =0$ for some $\alpha\in R$, we 
have to argue more carefully: Choose an open ball $B\subseteq\Omega$ with 
center $x_0$. Then $\sigma_\alpha x \in B$ for $x\in B$, and 
$\,\sigma_\alpha x - x = -\langle \alpha,x\rangle \alpha. $  Now Taylor's 
formula yields 
\[ f(\sigma_\alpha x) - f(x) = -\langle\alpha ,x\rangle\,\langle\nabla f(x),
\alpha\rangle + \frac{1}{2}\langle \alpha, x\rangle^2 \,
\alpha^t D^2f(\xi)\alpha\,,\]
with some $\xi$ on the line segment between $x$ and $\sigma_\alpha x$. It 
follows that for $x\in B$ with $\langle\alpha ,x\rangle \not=0$ we have 
$\,\delta_\alpha f(x) = \frac{1}{2}\alpha^t D^2u(\xi)\alpha.$ Passing to the 
limit $x\to x_0$ now leads to $\delta_\alpha f(x_0) = \frac{1}{2}\alpha^t D^2 
f(x_0) \alpha \leq 0$, which finishes the proof. 
\end{proof}

At this stage it is not much effort to gain a weak maximum principle for 
$k$-subharmonic functions on bounded, $G$-invariant subsets of $\b R^N$, 
which we want to include  here before passing over to the heat equation. Its 
range of validity is quite general, in contrast to the strong maximum 
principle in \cite{Du1}, which is restricted to $k$-harmonic polynomials on 
the unit ball. Our proof follows the classical one for the usual Laplacian, 
as it can be found e.g. in \cite{Jo}. 

\begin{theorem} Let $\Omega\subset \b R^N$ be open, bounded and 
$G$-invariant, and let $f\in C^2(\Omega)\cap C(\overline\Omega)$ be 
real-valued and $k$-subharmonic on $\Omega$. Then
\[ {\rm max}_{\,\overline \Omega} (f)\,=\, {\rm \max}_{\,\partial\Omega} 
(f)\,.\]
\end{theorem}

\begin{proof}
Fix  $\epsilon>0$ and put $g:= f + \epsilon |x|^2$. A short calculation 
gives $\Delta_k(|x|^2) = 2N+4\gamma >0.$ Hence $\Delta_k g >0$ on $\Omega$, 
and Lemma 4.1 shows that $g$ cannot achieve its maximum on $\overline\Omega$ 
at any $x_0\in\Omega.$ It follows that
\[ {\rm max}_{\,\overline\Omega}\,(f + \epsilon|x|^2)\,=\, 
{\rm max}_{\,\partial\Omega}\,(f +\epsilon |x|^2)\]
for each $\epsilon >0$. Consequently,  
\[ {\rm \max}_{\,\overline\Omega} \,(f)\, +\, \epsilon \,
{\rm min}_{\,\overline\Omega} |x|^2 \,\leq\,{\rm max}_{\,\partial\Omega} (f)
 \,+\,\epsilon\, {\rm max}_{\,\partial\Omega} |x|^2.\]
The assertion now follows with $\epsilon\to 0$.
\end{proof}

In this section we consider the generalized heat operator
\[ H_k \,:=\, \Delta_k -\partial_t\]
on function spaces $C^2(\Omega\times(0,T)), $ where $T>0$ 
 and $\Omega\subseteq\b R^N$ is open and $G$-invariant.  
Among the variety  of initial- and boundary value problems which may be 
posed for  $H_k$ in analogy to the corresponding classical problems, we 
here focus on the homogeneous Cauchy problem:

Find $u\in C^2(\b R^N\times(0,T))$ which is continuous on $\b R^N\times[0,T]$ 
and satisfies
\begin{equation}
\begin{cases}\label{(2.15)} 
H_k u = 0 & \text{on $\b R^N\times(0,T)$},\\
     u(.\,,0) = f & \text{$\in C_b(\b R^N).$}
   \end{cases}
\end{equation}

First of all, let us note  some basic solutions of the generalized heat 
equation $H_ku =0$.  Again we set $\,\gamma:= \sum_{\alpha\in R_+} k(\alpha)\,
\geq 0.$ 

\begin{lemma}
For parameters $a\geq 0$ and $b\in \b R\setminus \{0\}$, the function
\[ u(x,t) = \frac{1}{(a-bt)^{\gamma+N/2}} \,\,
{\rm exp}\!\left\{\frac{b|x|^2}{4(a-bt)}\right\}\]
solves $\,H_k u = 0\,$ on $\displaystyle \b R^N \times(-\infty,a/b)$ in case 
$b>0$, and on $\displaystyle \b R^N\times (a/b, \infty)$ in case $b<0$.
\end{lemma}

\begin{proof} The product rule (2.1) together with  $\sum_{i=1}^N T_i x_i = 
N+2\gamma\,$ shows that for each $\lambda>0$,
\[ \Delta_k\bigl(e^{\lambda|x|^2}\bigr) = \sum_{i=1}^N T_i 
\bigl(2\lambda x_i\, e^{\lambda|x|^2}\bigr) = 2\lambda\,
(N+2\gamma +2\lambda|x|^2)\, e^{\lambda|x|^2}.\]
From this the statement is  obtained readily by a short calculation.
\end{proof}

In particular, the function 
\[ F_k(x,t) = \frac{M_k}{t^{\,\gamma+N/2}} e^{-|x|^2/4t}\,,
\quad\text{with}\quad 
M_k = 4^{-\gamma-N/2}c_k\,,\]
 is a solution of the  heat equation $H_k u = 0$ on $\b R^N\times 
(0,\infty).$ 
It generalizes the fundamental solution for the classical heat equation 
$\,\Delta u -\partial_t u =0$, which is given by $\, F_0(x,t) = 
(4\pi t)^{-N/2} e^{-|x|^2/4t}\,.$  The normalization constant $M_k$ is 
choosen such that 
\[\int_{\b R^N} F_k(x,t)\, w_k(x) dx\,=\, 1 \quad \text{for all}\>\> t>0.\]

In order to solve the Cauchy problem \eqref{(2.15)}, it suggests itself to 
apply Fourier transform methods -- in our case, the Dunkl transform -- under 
suitable decay assumptions on the initial data $f$. In fact, in the classical 
case $k=0$ a bounded solution of \eqref{(2.15)} is obtained by convolving 
$f$ with the fundamental solution $F_0$, and its uniqueness is a consequence 
of a well-known maximum principle for the heat operator. It is not much 
effort to extend this maximum principle to the generalized heat operator 
$H_k$  in order to obtain uniqeness results; we shall do this in 
Prop. 4.12 and Theorem 4.13 at the end of this section. However, in our 
general situation it is not known whether there 
 exists a reasonable convolution structure on $\b R^N$ matching the action 
of the Dunkl transform $\mathcal D_k$, i.e. making it a homomorphism on  
suitable function spaces. In the one-dimensional case this is true: there 
is a $L^1$-convolution algebra associated with the reflection group $\b Z_2$ 
on $\b R$ and the multiplicity  parameter $k=\mu\geq 0$; this convolution  
enjoys many properties of a group convolution. It is studied in \cite{Roe} 
(see also \cite{RV} and  \cite{Ro}). 

In the $N$-dimensional case, we may introduce the notion of a generalized 
translation at least on the Schwartz space $\mathcal S(\b R^N)$ (and similar 
on $L^2(\b R^N, w_k(x)dx)$, as follows:
\begin{equation}\label{(4.75)}
 L_k^y f(x):= \frac{c_k^2}{4^{\gamma+N/2}} \int_{\b R^N} \mathcal D_k f(\xi) 
\,K(ix,\xi) K(iy, \xi)\, w_k(\xi) d\xi; \quad y\in \b R^N, \> 
f\in \mathcal S(\b R^N).
\end{equation}
Note that in case $k=0$, we simply have $\,L_0^y f(x) = f(x+y),$ while in 
the one-dimensional case, \eqref{(4.75)} matches the above-mentioned 
convolution structure on $\b R.$   Clearly, $L_k^y f(x) = L_k^x f(y)$; 
moreover, the inversion theorem for the Dunkl transform assures that 
$L_k^y f = f$ for $y=0$ and 
$\,\mathcal D_k(L_k^y f)(\xi) = K(iy,\xi) \mathcal D_k f(\xi).\,$
 From this it is not hard to see (by use of the bounds (2.7)) that 
$L_k^y f$ belongs to $\mathcal  S(\b R^N)$ again. 

Let us now consider the ``fundamental solution'' $F_k(.\,,t)$ for $t>0$. A 
short calculation, using Prop. 3.10 or Lemma 4.11 of \cite{deJ}, shows that
\begin{equation}\label{(4.80)}
 (\mathcal D_k F_k)(\xi,t)\,=\, e^{-t|\xi|^2}.
\end{equation}
By use of the reproducing formula (2.8) one therefore obtains from 
\eqref{(4.75)} the representation
\begin{equation}\label{(4.81)}
 L_k^{-y} F_k(x,t)\,=\, \frac{M_k}{t^{\gamma+N/2}} e^{-(|x|^2 + |y|^2)/4t}\,
K\Bigl(\frac{x}{\sqrt{2t}},\frac{y}{\sqrt{2t}}\Bigr).
\end{equation}

\begin{definition}
The generalized heat kernel $\Gamma_k$ is given by
\[ \Gamma_k(x,y,t):=\, \frac{M_k}{t^{\gamma +N/2}}\,e^{-(|x|^2 + |y|^2)/4t}\,
K\Bigl(\frac{x}{\sqrt{2t}},\frac{y}{\sqrt{2t}}\Bigr),\quad x,y\in \b R^N,\> 
t>0.\]
\end{definition}

\begin{lemma} The heat kernel $\Gamma_k$ has the following properties on 
$\b R^N\times\b R^N\times (0,\infty)$: \parskip=-2pt
\begin{enumerate}\itemsep=2pt
\item[\rm{(1)}] $\,\displaystyle  \Gamma_k(x,y,t) \,=\, 
\frac{c_k^2}{4^{\gamma+N/2}} \int_{\b R^N} e^{-t|\xi|^2}\, K(ix,\xi)\, 
K(-iy,\xi) \, w_k(\xi) d\xi\,.$
\item[\rm{(2)}] For fixed $y\in \b R^N$, the function $u(x,t):= 
\Gamma_k(x,y,t)$ solves the generalized heat equation $\, H_ku=0$ on 
$\b R^N\times (0,\infty)$.
\item[\rm{(3)}] $\,\displaystyle \int_{\b R^N} \Gamma_k(x,y,t)\,w_k(x) 
dx\,=\,1.$
\item[\rm{(4)}] $\,\displaystyle |\Gamma_k(x,y,t)|\,\leq\,
\frac{M_k}{t^{\gamma+N/2}}\,e^{-(|x|-|y|)^2/4t}.$  
\end{enumerate}
\end{lemma}

\begin{proof}
(1) is clear from the above derivation. For (2), remember that 
$\,\Delta_k^x \,K(ix, \xi) = -|\xi|^2 K(ix,\xi).$ Hence the assertion  
follows at once from representation (1) by taking the differentiations 
under the integral. This is justified by the decay properties of the 
integrand and its derivatives in question (use estimation (2.7) for the 
partial derivatives of $K(ix,\xi)$ with respect to $x$.)  
To obtain (3), we employ formula \eqref{(4.81)} as well as \eqref{(4.80)} 
and  write
\[ \int_{\b R^N} \Gamma_k(x,y,t)\,w_k(x)dx\,=\, 
\mathcal D_k\bigl( L_k^{-y} F_k\bigr)(0,t)\,=\, 
K(-iy,0)(\mathcal D_k F_k) (0,t)\,=\, 1.\]
Finally, (4) is a a consequence of the estimate (2.4) for $K$. 
\end{proof}

\begin{remark} In contrast to the classical case, it is not yet clear at 
this stage that the kernel $\Gamma_k$ is generally nonnegative. In fact, 
it is still an open conjecture  that the function $\,K(iy,.)\,$ is  
positive-definite on $\b R^N$ for each $y\in \b R^N$ (c.f. the remarks in 
\cite{deJ} and \cite{Du3}). This would imply a Bochner-type integral 
representation of $K(iy,.)$ and positivity of $K$ on $\b R^N\times\b R^N$ 
as an immediate consequence. In the one-dimensional case this conjecture is 
true, and the Bochner-type integral representation is explicitely known 
(see \cite{Ro} or \cite{Roe}). By  one-parameter semigroup techniques, it 
will however soon turn out that $K$ is at least positive on 
$\b R^N\times \b R^N$. 
\end{remark}

\begin{definition}
For $f\in C_b(\b R^N)$ and $t\geq 0$ set 
\begin{equation}\label{(3.16)}
 H(t)f(x):= \begin{cases}
    \displaystyle\int_{\b R^N} \Gamma_k(x,y,t) f(y)\,w_k(y)dy & 
    \text{if $\,\,t>0$},\\
     f(x) & \text{if $\,\,t=0$}
   \end{cases}
\end{equation}
\end{definition}

Part (4) of Lemma 4.5 assures that for each $t\geq 0,\>\>  H(t)f$ is 
well-defined and continuous on $\b R^N$. It provides a natural candidate 
for the solution to our Cauchy problem. Indeed, when restricting to initial 
data from the Schwartz space $\mathcal S(\b R^N)$, we easily obtain the 
following:

\begin{theorem}
Suppose that $f\in\mathcal S(\b R^N).$ Then $\,u(x,t):= H(t)f(x), \>\, (x,t) 
\in \b R^N\times [0,\infty),$ 
 solves the Cauchy-problem \eqref{(2.15)} for each $T>0$. 
Morover, it has the following properties: \parskip=-1pt
\begin{enumerate}\itemsep=0pt
\item[\rm{(i)}] $\,H(t)f \in \mathcal S(\b R^N)$ for each $t>0$. 
\item[\rm{(ii)}] $\, H(t+s)f\,=\, H(t) H(s) f\>\>$ for all $s,t\geq 0$.
\item[\rm{(iii)}] $\, \|H(t)f - f\|_{\infty,\b R^N} \, \to 0\> $ with 
$t\to 0$. 
\end{enumerate}
\end{theorem}

\begin{proof}
Using formula (1) of Lemma 4.5  and Fubini's theorem, we can write
\begin{align}\label{(3.14)}
u(x,t) = H(t)f(x)\,=\,&\frac{c_k^2}{4^{\gamma+N/2}}\int_{\b R^N}
\int_{\b R^N} K(ix,\xi) K(-iy,\xi)\, e^{-t|\xi|^2} f(y)\,w_k(\xi)w_k(y)\,
d\xi dy \notag\\
=\,& \frac{c_k^2}{4^{\gamma+N/2}} \int_{\b R^N} e^{-t|\xi|^2}\mathcal D_k 
f(\xi) K(ix,\xi)\,  w_k(\xi) d\xi
\end{align}
for all $t>0$. (Remember that $\mathcal S(\b R^N)$ is invariant under the 
Dunkl transform). This makes clear that (i) is satisfied. As before, it is 
seen  that differentiations may be taken under the integral in 
\eqref{(3.14)}, and that  $H_k u = 0$ on $\b R^N\times (0,\infty)$. 
 Moreover, in view of the inversion theorem for the Dunkl transform, 
\eqref{(3.14)} holds for $t=0$ as well.  Using (2.5), we thus obtain  the 
estimation
\[ \|H(t)f - f\|_{\infty,\b R^N}\,\leq\, \sqrt{|G|} 
\frac{c_k^2}{4^{\gamma+N/2}}\int_{\b R^N} |\mathcal D_k f(\xi)| 
\bigl(1- e^{-t|\xi|^2}\bigr) w_k(\xi)d\xi,\]
and this integral tends to $0$ with $t\to 0.$ This yields (iii). In 
particular, it follows that $u$ is continuous on $\b R^N\times [0,\infty)$. 
To prove (ii), note that $\,\mathcal D_k\bigl(H(t)f\bigr)(\xi)\,=\, 
e^{-t|\xi|^2}\mathcal D_k f(\xi)\,.$ Therefore
\[\mathcal D_k\bigl(H(t+s)f\bigr)(\xi)\,=\, 
e^{-t|\xi|^2}\mathcal D_k\bigl(H(s)f\bigr)f(\xi)\,=\, 
\mathcal D_k\bigl( H(t)H(s)f\bigr)(\xi).\]
The statement now follows from the injectivity of the Dunkl transform on 
$\mathcal S(\b R^N)$. 
\end{proof}

We are now going to show that indeed, the linear operators $H(t)$ on 
$\mathcal S(\b R^N)$ extend to a  positive contraction  semigroup on the 
Banach space $C_0(\b R^N),$ equipped with its uniform norm $ \|.\|_\infty$. 
To this end, we consider the generalized Laplacian $\Delta_k$ as a densely 
defined linear operator on $C_0(\b R^N)$ with domain $\mathcal S(\b R^N)$.

\begin{theorem}\parskip=-1pt
\begin{enumerate}\itemsep=0pt
\item[\rm{(1)}] The operator $\Delta_k$ on $C_0(\b R^N)$  is closable, and 
its closure $\overline\Delta_k$ generates a positive, strongly continuous 
contraction semigroup $\{T(t),t\geq 0\}$ on $C_0(\b R^N)$.
\item[\rm{(2)}] The action of $T(t)$ on $\mathcal S(\b R^N)$ is given by 
$\, T(t)f = H(t)f$.
\end{enumerate}
\end{theorem}

\begin{proof} (1) We apply a variant of the Lumer-Phillips theorem, which 
characterizes  generators of positive one-parameter contraction semigroups 
(see e.g. \cite{Ar}, Cor. 1.3). It requires two properties: \parskip=-1pt
\begin{enumerate}\itemsep=0pt
\item[{\rm(i)}] The operator $\Delta_k$ satisfies the following 
``dispersivity condition'':\enskip  
 Suppose that $f\in \mathcal S(\b R^N)$ is real-valued with 
$\,{\rm max}\{f(x): x\in \b R^N\} = f(x_0).$ Then $\Delta_k f(x_0) \leq 0.$

\item[{\rm(ii)}] The range of $\lambda I  -\Delta_k$ is dense in 
$C_0(\b R^N)$ for some $\lambda >0.$ 
\end{enumerate}
Property (i) is  an immediate consequence of Lemma 4.1. Condition (ii) is 
also satisfied, because $\lambda I  -\Delta_k$ maps $\mathcal S(\b R^N)$ 
onto itself for each $\lambda >0$;  this follows from the fact that the 
Dunkl transform  is a  homeomorphism of $\mathcal S(\b R^N)$ and 
$\,\mathcal D_k\bigl((\lambda I  -\Delta_k)f\bigr)(\xi) = 
(\lambda +|\xi|^2)\mathcal D_k f(\xi)$. 
The assertion now follows by the above-mentioned theorem.

(2) It is known from semigroup theory that for every 
$f\in \mathcal S(\b R^N),$ the function 
$t\mapsto T(t)f$ is the unique solution of the abstract Cauchy problem
\begin{equation}
\begin{cases}\label{(4.70)} 
\displaystyle\frac{d}{dt}u(t)= \overline\Delta_k u(t) & \text{for $t>0$},\\
     u(0) = f & \text{}
   \end{cases}
\end{equation}
within the class of all (strongly) continuously differentiable functions $u$ 
on $[0,\infty)$ with values in the Banach space $(C_0(\b R^N), \|.\|_\infty).$
  By property (i) of Theorem 4.6 we have $H(t)f\in C_0(\b R^N)$ for $f\in 
\mathcal S(\b R^N)$. Moreover, from representation \eqref{(3.14)} of $H(t)f$ 
it is  readily seen that $\,t\mapsto H(t)f$ is continuously differentiable 
on $[0,\infty)$ and solves \eqref{(4.70)}.
This finishes the proof.
\end{proof}

\begin{corollary}
The  heat kernel $\Gamma_k$ is strictly positive on 
$\b R^N\times\b R^N\times(0,\infty).$ In particular, the generalized 
exponential kernel K satisfies 
\[ K(x,y) > 0 \quad \text{for all}\>\> x,y\in \b R^N.\]
\end{corollary}

\begin{proof}
For any  initial distribution $f\in \mathcal S(\b R^N)$ with $f\geq 0$ the 
last theorem  implies that 
\[\int_{\b R^N} \Gamma_k(x,y,t) f(y)\,w_k(y)dy \,= \, T(t)f(x)\,\geq\, 0 
\quad \text{for all}\>\> (x,t)\in \b R^N\times[0,\infty).\]
As $\,y\mapsto \Gamma_k(x,y,t)$ is continuous on $\b R^N$ for each fixed 
$x\in \b R^N$ and $t>0$, it follows that $\Gamma_k(x,y,t) \geq 0$ for all 
$x,y\in \b R^N$ and $t>0.$ Hence $K$ is nonnegative as well. 
Now recall again the reproducing identity (2.8), which says that
\[ e^{(|x|^2 +|y|^2)} K(2x,y) \,=\,  c_k \int_{\b R^N} K(x,2z) K(y,2z) 
\,w_k(z) \,e^{-|z|^2} dz\]
for all $x,y\in \b R^N.$ The integrand on the right side is continuous, 
non-negative and not identically zero (because $K(x,0) K(y,0) =1$). Therefore 
the integral itself must be strictly positive. 
\end{proof}

\begin{corollary}
The semigroup $\{T(t)\}$ on $C_0(\b R^N)$ is given explicitely by 
\[ T(t)f\,=\, H(t)f, \>\> f\in C_0(\b R^N).\]
\end{corollary}

\begin{proof} This is clear from part (2) of  Theorem 4.8 and the previous 
corollary, which implies that the operators  $H(t)$ are continuous -- even 
contractive --  on $C_0(\b R^N).$ 
\end{proof}

\begin{remark} The generalized Laplacian also leads to a contraction 
semigroup on the Hilbert space $\mathcal H:= L^2(\b R^N, w_k(x)dx);$ this 
generalizes the results of \cite{Ro} for the  one-dimensional case. In fact, 
let $M$ denote the multiplication operator on $\mathcal H$ defined by 
$\, Mf(x) = -|x|^2 f(x)\,$ and with domain $\, D(M) = \{f\in \mathcal H: 
|x|^2 f(x)\in \mathcal H\}$.
$M$ is self-adjoint and generates the strongly continuous contraction 
semigroup $\,M(t)f(x) = e^{-t|x|^2}f(x) \> (t\geq 0)$ on $\mathcal H$.  For 
$f\in \mathcal S(\b R^N)$, we have the identity $\,\mathcal D_k(\Delta_kf) = 
M(\mathcal D_k f).$ As $\mathcal S(\b R^N)$ is dense in $D(M),$ this shows 
that $\Delta_k$ has a self-adjoint extension $\widetilde\Delta_k$ on 
$\mathcal H$, namely $\,\widetilde\Delta_k = 
\mathcal D_k^{-1} M \mathcal D_k$, where here $\mathcal D_k$ denotes the 
Plancherel-extension of the Dunkl transform to $\mathcal H$. The domain of 
$\widetilde\Delta_k$ is the Sobolev-type space $\, D(\widetilde \Delta_k) =
\{f\in \mathcal H: |\xi|^2\,\mathcal D_k f(\xi)\in \mathcal H\}.$ Being  
unitarily equivalent with $M,$ the opertor $\widetilde\Delta_k$  also 
generates a strongly continuous contraction semigroup $\{\widetilde T(t)\}$ 
on $\mathcal H$ which is unitarily equivalent with $\{M(t)\};$ it is given by
\[\widetilde T(t)f(x)\,=\,\int_{\b R^N} e^{-t|\xi|^2} \mathcal D_k f(\xi)\, 
K(ix,\xi)\,w_k(\xi)d\xi.\]
\end{remark}

The knowledge that $\Gamma_k$ is nonnegative allows also to solve the Cauchy 
problem  \eqref{(2.15)} in its general setting:

\begin{theorem}
Let $f\in C_b(\b R^N)$. Then $\,u(x,t):= H(t)f(x)$ is bounded on 
$\b R^N\times [0,\infty)$ and solves the Cauchy problem \eqref{(2.15)} for 
each $T>0.$ 
\end{theorem}

\begin{proof}
In order to see that $u$ is twice continuously differentiable on 
$\b R^N\times (0,\infty)$ with $\,H_k u=0$, we only have to make sure that 
the necessary differentiations of $u$ may be taken under the integral in 
\eqref{(3.16)}. One has to use again the estimations (2.7) for the partial 
derivatives of $K$; these  provide sufficient decay properties of the 
derivatives of $\Gamma_k$, allowing the necessary differentiations of $u$ 
under the integral by use of the dominated convergence theorem. 
Boundedness of $u$ is clear from the positivity and normalization 
(Lemma 4.5(3))  of $\Gamma_k$; in fact, 
$\, |u(x,t)|\,\leq\, \|f\|_{\infty, \b R^N}\,$ on $  \b R^N\times 
[0,\infty).$

Finally, we have to show that $\,H(t)f(x)\to f(\xi)\,$ with $x\to\xi$ and 
$t\to 0$. We start by the usual method: For fixed $\epsilon>0$, choose  
$\delta >0$ such that $|f(y)- f(\xi)| <\epsilon $ for $|y-\xi| <2\delta$ 
and let $M:= \|f\|_{\infty,\b R^N}.$ Keeping  in mind the positivity and 
normalization  of $\Gamma_k$, we obtain for $|x-\xi| <\delta$ the estimation
\begin{align}
 |H(t)f(x)- &f(\xi)|\,\leq\,\Big|\int_{\b R^N} 
\Gamma_k(x,y,t)\big( f(y)-f(\xi)\big) w_k(y)dy\Big| \notag\\
 \leq\,& \int_{|y-x|<\delta} \Gamma_k(x,y,t) |f(y)- f(\xi)|w_k(y)dy\,+\, 
\int_{|y-x|>\delta} \Gamma_k(x,y,t)|f(y)- f(\xi)| w_k(y)dy 
\notag\\
 <\,& \epsilon +\,2M \int_{|y-x|>\delta} \Gamma_k(x,y,t) w_k(y)dy. \notag
\end{align}
It thus remains to show that for each $\delta>0$, 
\[ \text{lim}_{(x,t)\to(\xi,0)} \int_{|y-x|>\delta} \Gamma_k(x,y,t) 
w_k(y)dy\,=\,0.\]
For abbreviation put
\[ I(x,t):= \int_{|y-x|\leq\delta} \Gamma_k(x,y,t) w_k(y) dy.\]
As $I(x,t)\leq 1$, it suffices to prove that 
$\,\liminf_{(x,t)\to(\xi,0)} I(x,t)\,\geq\,1.\,$
For this, choose some positive constant $\delta^\prime <\delta$ and 
$h\in\mathcal S(\b R^N)$ with $\,0\leq h\leq 1, \> h(\xi)=1$ and such 
that $h(y) =0$ for all $y$ with $|y-\xi|>\delta-\delta^\prime.$ Then for 
each $x$ with $|x-\xi|<\delta^\prime$ the support of $h$ is contained in 
$\{y\in \b R^N: |y-x|\leq\delta\}$; therefore
\[ \int_{\b R^N} h(y)\Gamma_k(x,y,t) w_k(y)dy\,\leq\, I(x,t)\]
for all $(x,t)$ with $|x-\xi|<\delta^\prime.$ But according to Theorem 4.7 
we have
\[ \text{lim}_{(x,t)\to(\xi,0)} \int_{\b R^N} h(y)\Gamma_k(x,y,t)\,w_k(y)dy 
\,=\, h(\xi)\,=\,1.\]
This finishes the proof. 
\end{proof}

It is still open whether our solution of the Cauchy problem \eqref{(2.15)} 
is unique within an appropriate class of functions. As in the classical case, 
this follows from an maximum principle for the generalized heat operator  
on $\b R^N\times (0,\infty)$. The first step is the following weak maximum 
principle for $H_k$  on bounded domains. It is proved by  a similar method 
as used in Theorem 4.2. By virtue of Lemma 4.1, this proof  is litterally 
the same as the standard proof in the classical case (see  e.g. \cite{Jo}) 
and therefore omitted here.

\begin{proposition}  Suppose that $\Omega\subset\b R^N$ is open, bounded 
and $G$-invariant. For $T>0$ set
\[\Omega_T := \Omega\times(0,\infty)\quad \text{and}\quad\,   
\partial_*\Omega_T := \left\{(x,t)\in \partial \Omega_T: \,t=0 \>\, 
\text{or} \>\>x\in\partial\Omega\,\right\}.\]
Assume further that $u\in C^2(\Omega_T)\cap C(\overline\Omega_T)$ 
satisfies $\,H_k u \geq 0$ in $\Omega_T$. Then
\[ {\rm max}_{\,\overline\Omega_T} (u)\,=\, {\rm max}_{\,\partial_*\Omega_T} 
(u)\,.\]
\end{proposition}

Under a suitable growth condition on the solution, this maximum principle 
may be extended to the case where $\Omega=\b R^N$. The proof is adapted from 
the one in \cite{dBe} for the classical case. 

\begin{theorem} {\rm (Weak maximum principle for $H_k$ on $\b R^N$.)} 
\enskip\enskip Let $S_T := \b R^N\times (0,T)$ and suppose that 
$\,u\in C^2(S_T)\cap C(\overline S_T)$ satisfies 
\[ \begin{cases}
     H_k u \geq 0 & \text{in $S_T$},\\
     u(.\,,0) = f\,, & 
   \end{cases}\]
where $f\in C_b(\b R^N)$ is real-valued. Assume further that 
there exist positive constants $C, \lambda, r$ such that
\[ u(x,t) \,\leq\, C\cdot e^{\lambda|x|^2} \quad\text{for all} \>\,
(x,t)\in S_T \quad\text{with}\>\, |x|>r.\]
 Then
\(\displaystyle \quad\, {\rm max}_{\,\overline S_T} (u)\,\leq\, 
{\rm sup}_{\b R^N} (f).\)
\end{theorem}

\begin{proof}
Let us first assume  that $8\lambda T < 1.$ For fixed $\epsilon >0$ set
\[ v(x,t) := u(x,t) - \epsilon\cdot\frac{1}{(2T-t)^{\,\gamma+N/2}}\,
{\rm exp}\left\{\frac{|x|^2}{4(2T-t)}\right\}, \quad (x,t) 
\in \b R^N\times [0,2T).\]
By Lemma 4.3, $v$ satisfies
\(\, H_k v\,=\, H_k u \,\geq 0\,\) in $S_T$. 
Now fix some constant $\rho>r$ and consider the bounded cylinder 
$\Omega_T = \Omega \times(0,T)$ with $\Omega = 
\{x\in \b R^N:\,|x|<\rho\}.$ Setting $\, M:= \sup_{\b R^N}(f)$, we have
\(\, v(x,0) < u(x,0) \leq M \,\) for $x\in\overline \Omega.$ 
Moreover, for $\,|x|=\rho$ and $t\in (0,T]$ 
\[ v(x,t) \,\leq\, C e^{\lambda \rho^2} - 
\epsilon\cdot\frac{1}{(2T)^{\,\gamma+N/2}}\,e^{\,\rho^2/8T}.\]
As $\,\lambda< (8T)^{-1}$, we see that $\,v(x,t)\leq M\,$ on 
$\partial_*\Omega_T$, provided that $\rho$ is large enough. Then by 
Prop. 4.12  we also have $\,v(x,t)\leq M\,$ on $\overline\Omega_T\,.$ As 
$\rho>r$ was arbitrary, it follows that $\,v(x,t)\leq M\,$ on 
$\overline S_T$. As $\epsilon>0$ was arbitrary as well, this implies that 
$\,u(x,t)\leq M\,$ on $\overline S_T$. If $8\lambda T \geq 1$, we may 
subdivide $S_T$ into finitely many adjacent open strips of width less 
than $1/8\lambda$ and apply the above conclusion repeatedly. 
\end{proof}
 
\begin{corollary}
The solution   of the Cauchy problem (4.1) according to Theorem 4.11 is 
unique within the class of functions $ u\in C^2(S_T)\cap C(\overline S_T)$ 
which satisfy the following exponential growth condition:\enskip
There exist positive constants $\,C,\lambda, r\,$ such that 
\[ |u(x,t)|\,\leq\, C\cdot e^{\lambda |x|^2} \quad\text{for all}\>\> 
(x,t)\in S_T \>\>\text{with}\>\>|x|>r.\] 
\end{corollary}

\end{document}